\documentclass[11pt]{amsart}
\usepackage{amsmath}
\usepackage{amssymb}
\usepackage{amsthm}
\usepackage{mathrsfs}
\usepackage{enumerate}
\usepackage[english]{babel}
\usepackage{hyperref}
\usepackage{xypic}

\newtheorem{thm}{Theorem}[section]
\newtheorem{lem}[thm]{Lemma}
\newtheorem{cor}[thm]{Corollary}

\theoremstyle{definition}

\newtheorem{rem}[thm]{Remark}
\newtheorem{defn}[thm]{Definition}

\newtheorem{ex}[thm]{Example}

\newtheorem*{acknowledgments*}{Acknowledgments}

\numberwithin{equation}{section}

\theoremstyle{remark}
\newtheorem{notation}[thm]{Notation}

\newcommand\bk[1]{\langle #1 \rangle}

\mathchardef\ordinarycolon\mathcode`\: 
\def\vcentcolon{\mathrel{\mathop\ordinarycolon}} 
\providecommand*\coloneqq{\mathrel{\vcentcolon\mkern-1.2mu}=}

\def\Z{{\mathbb Z}} 
\def\R{{\mathbb R}} 
\def\C{{\mathbb C}} 
\def\N{{\mathbb N}}

\DeclareMathOperator\Imaginary{Im}
\DeclareMathOperator\Real{Re}

\def\A{{\mathcal A}}
\def\H{{\mathcal H}}
\def\D{\mathscr D}
\def\AHD{{(\A, \H, \D)}}

\def\Spec{{\mathrm{Spec}}}

\def\alg{\mathrm{alg}}

\def\dom{{\mathrm{dom}}}

\def\PDO{\Psi}
\def\BPDO{\mathcal B}
\def\DO{\mathcal D}

\def\CL{\Psi_{\mathrm{cl}}}
\newcommand\W[1]{W^{#1}}
\def\Linear{\mathcal L}
\def\Compact{\mathcal K}
\newcommand\Op[1]{\mathrm{Op}^{#1}}
\def\ord{\mathrm{order}}

\newcommand\out[1]{\!}

\def\End{{\mathrm{End}}}
\def\tensor{\widehat{\otimes }}

\newcommand\half[1]{\frac{#1}{2}}

\begin{document}

\title[$\Psi$DO's and Regularity]{Pseudo-differential Operators and Regularity of Spectral Triples}
\author{Otgonbayar Uuye}
\date{July 30, 2010}
\address{
Department of Mathematical Sciences\\
University of Copenhagen\\
Universitetsparken 5\\
DK-2100 Copenhagen E\\
Denmark}
\email{otogo@math.ku.dk}
\urladdr{http://www.math.ku.dk/~otogo}
\maketitle

\begin{abstract}
We introduce a notion of an algebra of generalized pseudo-differential operators and prove that a spectral triple is regular if and only if it admits an algebra of generalized pseudo-differential operators. We also provide a self-contained proof of the fact that the product of regular spectral triples is regular.
\end{abstract}

\setcounter{section}{-1}
\section{Introduction}

In their work on the local index theorem \cite{CM:1995}, Connes and Moscovici introduced the notion of a regular spectral triple. Later, Higson showed that a spectral triple is regular if and only if a certain algebra constructed from the spectral triple is what he calls an algebra of generalized differential operators \cite{H:2004a, H:2006}. See Theorem~\ref{Higson Regularity} for the precise statement.

Motivated by their work we introduce a notion of an algebra of generalized pseudo-differential operators; it clarifies the role of the regularity condition and serves as a convenient framework to study index theory via complex powers. I should note that most of the main ideas in this paper can be traced back to \cite[Appendix B]{CM:1995} or \cite[Subsection 4.5]{H:2006}, but for the convenience of the reader we tried to remain self-contained.

Now we describe the content of the paper. After some preliminaries in Section \ref{Prelim}, we recall Higson's notion of an algebra of generalized differential operators in Section~\ref{Algebra of Generalized Differential Operators}. In Section \ref{Algebra of Generalized Pseudodifferential Operators}, we develop the notion of an algebra of generalized pseudo-differential operators. As usual, there are two versions -- operators of all orders and operators of order at most zero. Extending Higson's result, we show that a spectral triple is regular if and only if it admits an algebra of generalized pseudo-differential operators (see Theorem~\ref{thm regularity}).

Finally, we show that the product of regular spectral triples is again regular (see Theorem~\ref{thm: product of regular}). This is a folklore and it is implicitly contained in or follows from results in \cite{CM:1995, H:2004a, H:2006}. But no direct reference seems to exist in the literature.

\begin{acknowledgments*} This paper is essentially contained in my thesis. I would like to thank my advisor Nigel Higson for his constant encouragement and guidance. I would also like to thank the NCG group at the University of Copenhagen for support.
\end{acknowledgments*}
  
\section{Preliminaries}\label{Prelim}

\subsection{Spectral Triples}

The following definition of a regular spectral triple is due to Connes and Moscovici \cite{CM:1995}.
\begin{defn}\label{defn spectral triple} A {\em spectral triple} $\AHD$  consists of an associative algebra $\A$, a graded Hilbert space $\H$ equipped with an even representation of $\A$ and a densely defined, self-adjoint, odd operator $\D$ such that 
\begin{enumerate}
\item $\A$ is contained in the domain of the derivation $[\D, -]$
, that is, if $\dom(\D)$ is the domain of $\D$,  then any $a \in \A$ satisfies $a \cdot \dom(\D) \subseteq \dom(\D)$ and the commutator $[\D, a]: \dom(\D) \to \H$ extends by continuity to a bounded operator on $\H$ and
\item\label{compact resolvent} 
the operators $a \cdot (\D \pm i)^{-1}$ and $(\D \pm i)^{-1} \cdot a$ are compact for any $a \in \A$.
\end{enumerate}

We say that a spectral triple $\AHD$ is {\em regular}, if
\begin{enumerate}
\item[(3)]\label{regularity CM} the space $\A + [\D, \A]$ is contained in the {\em smooth domain} of the derivation $\delta(T) \coloneqq [|\D|, T]$ on $\Linear(\H)$.
\end{enumerate}

\end{defn}

We remind the reader that the smooth domain of $\delta$ is defined as follows. 
The {\em domain}\index{derivation!domain of} $\dom(\delta)$ of $\delta = [|\D|, -]$ is the set of $T \in \Linear(\H)$ such that
	\begin{enumerate}[(i)] 
	\item $T \cdot \dom(|\D|) \subseteq \dom(|\D|)$ and
	\item $[|\D|, T]$ extends to a bounded operator on $\H$.
	\end{enumerate}
For $k \ge 2$, we define $\dom(\delta^{k})$ inductively as 
	\begin{equation}\dom(\delta^{k}) \coloneqq \{b \in \dom(\delta) \mid \delta(b) \in \dom(\delta^{k-1})\}.\end{equation} 
The smooth domain of $\delta$ is $\dom^{\infty}(\delta) \coloneqq \bigcap_{k=1}^{\infty}\dom(\delta^{k}).$

For simplicity, we do not consider any additional structure on $\A$; even though in most natural examples $\A$ have a topology or a norm, it will not play any role in our analysis.

Obviously, the following is the basic example that should be mentioned in any paper on spectral triples.
\begin{ex}\label{comm example} 
Let $M$ be a complete Riemannian manifold of even dimension and let $S \to M$ be a complex spinor bundle.  Let $\H \coloneqq L^{2}(M, S)$ denote the graded Hilbert space of $L^{2}$-sections of $S$. Let $\D$ be a Dirac-type operator acting on $C^{\infty}_{c}(M, S)$. Then standard $\Psi$DO theory implies that $\D$ with domain $C^{\infty}_{c}(M, S) \subset L^{2}(M, S)$ is essentially self-adjoint and 
	\begin{equation}
	(C^{\infty}_{c}(M), L^{2}(M, S), \bar\D),
	\end{equation}
is a {\em regular} spectral triple, where $\bar\D$ is the closure of $\D$. See \cite{HR:2000a}.
\end{ex}


The material in the following two subsections are elementary and well-known. We include it for reference and the convenience of the reader.

\subsection{Sobolev Spaces}\label{sec: Sobolev spaces}

Let $\Delta$ be an {\em invertible, positive, self-adjoint operator} on a Hilbert space $\H$. 

\begin{defn} 
The {\em $\Delta$-Sobolev space}\index{Sobolev space} of order $s \in \R$, denoted 
	\begin{equation}
	\W{s} = \W{s}(\Delta) = \W{s}(\H, \Delta),
	\end{equation} 
is the Hilbert space completion of $\dom(\Delta^{\half{s}})$ with respect to the inner product given by
	\begin{equation}\bk{\xi, \eta}_{\W{s}} \coloneqq \bk{\Delta^{\half{s}}\xi, \Delta^{\half{s}}\eta}_{\H}\end{equation}
for $\xi$, $\eta \in \dom(\Delta^{\half{s}})$.
\end{defn}
Note that the invertibility hypothesis guaranties that $||\cdot||_{\W{s}}$ is nondegenerate. The following is well-known.


\begin{lem}\label{lem: domain of power} For $s \ge t$, we have a continuous inclusion
	\begin{equation}
	\W{s} \subseteq \W{t}.
	\end{equation}
Moreover, for $s \ge 0$, $\dom(\Delta^{\half{s}})$ is complete and thus $\W{s} = \dom(\Delta^{\half{s}})$. \qed
\end{lem}  



\begin{ex} If $\Delta$ is {\em bounded}, then nothing interesting happens: for any $s \in \R$, $\W{s} = \H$ and the inner products are different but equivalent. 
\end{ex}


\begin{ex} Suppose that $\Delta$ has compact resolvents. 
 Then there exists a complete orthonormal basis $\{\xi_{n}\}$ of $\H$ consisting of eigenvectors: $\Delta\xi_{n} = \lambda_{n}\xi_{n}$ with $0 < \lambda_{1} \le \lambda_{2} \le \dots \to \infty.$
Hence the Sobolev spaces can be identified as the weighted $l^{2}$-space:
	\begin{equation}\W{s} = \{(a_{n})_{n=1}^{\infty} \mid \sum_{n} \lambda_{n}^{s}|a_{n}|^{2} < \infty\}.\end{equation}
\end{ex}

\begin{ex}[Classical Sobolev Spaces]\label{ex: classical Sobolev spaces} Let $M$ be a closed manifold and let $\Delta_{0}$ be a {\em strictly positive} order-two elliptic partial differential operator acting on the smooth functions on $M$. Let $\Delta = \bar\Delta_{0}$ be its closure. Then $\Delta$ is invertible and self-adjoint  and has compact resolvents (cf. \cite{HR:2000a}). It follows from the basic estimate that as a Hilbert space $\W{s}$ is equivalent to the ``standard'' $s$-Sobolev space (cf. \cite[Proposition I.7.3]{Shu:2001}).

\end{ex}

\begin{defn}Let $\Delta$ be an invertible, positive, self-adjoint operator. The space of {\em $\Delta$-smooth vectors}\index{smooth vectors} is
	\begin{equation}\W{\infty} \coloneqq \bigcap_{s \in \R} \W{s} = \bigcap_{n = 0}^{\infty} \W{2n} = \bigcap_{n=0}^{\infty} \dom(\Delta^{n}) \subseteq \H.\end{equation}
\end{defn}
The space $\W{\infty}$ contains $1_{[0, M]}(\Delta)\H$ for any $M > 0$, and thus $\W{\infty}$ is dense in $\dom(\Delta^{z})$ for any $z \in \C$. It follows that $\W{\infty}$ is dense in $\W{s}$ for any $s\in \R$.

\begin{lem}\label{lem: order of Delta^{z}}Let $z\in \C$. For any $s \in \R$, the operator $\Delta^{z}_{|\W{\infty}}: \W{\infty} \to \H$ extends\footnote{Note that if $s > 0$ the range actually shrinks to $\W{s} \subset \H$.} to an {\em isometry} 
	\begin{equation}\Delta^{z}: \W{s + 2\Real(z)} \to \W{s}.\end{equation}
In particular, $\Delta^{z}(\W{\infty}) = \W{\infty}$.
\end{lem}
\begin{proof} First note that $\Delta^{i\Imaginary(z)}$ is a unitary operator on $\H$. Therefore, for any $\xi \in \W{\infty}$,
	\begin{equation}||\Delta^{z}\xi||_{\W{s}} = ||\Delta^{\Real(z)+\half{s}}\xi||_{\H} = ||\xi||_{\W{s + 2\Real(z)}}.\end{equation}
\end{proof}

\begin{lem} The space $\W{\infty} \subset \H$ is a common core for the operators $\Delta^{z}$, $z \in \C$, {i.e.}\ $\Delta^{z}$ is essentially self-adjoint on $\W{\infty} \subset \H$.
\end{lem}
\begin{proof} Any vector $\xi \in 1_{[0, M]}(\Delta)\H \subset \W{\infty}$ is an analytical vector for $\Delta^{z}$. Indeed, for any $\xi \in 1_{[0, M]}(\Delta)\H$,
	\begin{equation}||(\Delta^{z})^{n}\xi|| \le M^{|\Real(z)n|}||\xi||\end{equation}
and consequently 
	\begin{equation}\sum_{n=0}^{\infty}\frac{||(\Delta^{z})^{n}\xi||}{n!}t^{n} < \infty\end{equation}
for any $t > 0$. Since $\Delta^{z}$ preserves $\W{\infty}$, all the analytical vectors for $\Delta^{z}$ are also analytical for the restriction $\Delta^{z}_{|\W{\infty}}$ and applying Nelson's theorem \cite[Theorem X.39]{RS:1975} to $\Delta^{z}_{|\W{\infty}}$, we see that $\Delta^{z}$ is essentially self-adjoint on $\W{\infty}$.
\end{proof}

Another convenient way to express the complex powers $\Delta^{z}$ is using the  Cauchy integral formula: for $\Real(z) < 0$,
	\begin{equation}\Delta^{z} = \frac{1}{2\pi i}\int \lambda^{z}(\lambda - \Delta)^{-1}d\lambda,\end{equation}
where the integral is a contour integral along a downwards pointing vertical line in $\C$ which separates $0$ from $\Spec(\Delta)$. Let $\Delta \ge c > 0$ and let $s \in \R$. Then it follows from the spectral theorem that for any $\lambda \notin \Spec(\Delta)$, the resolvent $(\lambda - \Delta)^{-1}$ is a bounded operator on the Sobolev spaces $\W{s}$ with norm at most $\left((\Real(\lambda) - c)^{2} + \Imaginary(\lambda)^{2}\right)^{-\half{1}}$. Hence for any $\Real(z) < 0$, the integral converges to a bounded operator on $\W{s}$.

More generally, we have the following.
\begin{lem}\label{lem: Cauchy integral formula for complex power} For $k \in \N$ and $\Real(z) < k$,
	\begin{equation}{z \choose k}\Delta^{z - k} = \frac{1}{2\pi i}\int \lambda^{z}(\lambda -\Delta)^{-k - 1}d\lambda\quad \text{in}\quad \Linear(\W{s}),\; s \in \R,\end{equation}
where the integral is a contour integral along a downwards pointing vertical line in $\C$ which separates $0$ from $\Spec(\Delta)$. 
\qed
\end{lem}

\subsection{Operators of Finite Analytic Order}
We consider various class of linear operators on $\W{\infty}$. The algebra of all linear operators $\W{\infty} \to \W{\infty}$ is denoted $\End(\W{\infty})$. 

\begin{ex}\label{ex: P gives an element of EndW} 
If an (unbounded) operator $P$
has domain $\dom(P) \supseteq \W{\infty}$ and preserves $\W{\infty}$, {\em i.e.}\ $P(\W{\infty}) \subseteq \W{\infty}$, then the restriction $P_{|\W{\infty}}$ gives an element of $\End(\W{\infty})$. We often write, simply, $P$ for $P_{|\W{\infty}}$.
\end{ex}

\begin{defn}\label{defn: analytic order} We say that a linear operator $\W{\infty} \to \W{\infty}$ has {\em analytic order at most $t \in \R$}\index{analytic order} if it extends by continuity to a {\em bounded} linear operator $\W{s + t} \to \W{s}$ for every $s \in \R$.
We write 
	\begin{equation}\Op{t} = \Op{t}(\Delta) = \Op{t}(\H, \Delta)\end{equation} 
for the class of operators of analytic order at most $t$ and define 
	\begin{equation}\Op{} = \Op{\infty} \coloneqq \bigcup_{t} \Op{t} \text{\: and \:} \Op{-\infty} \coloneqq \bigcap_{t} \Op{t}.\end{equation} 
\end{defn}

\begin{lem} Operators with finite analytic order form a filtered algebra:
	\begin{enumerate}[(a)]
	\item $\Op{s} \subseteq \Op{t}$ for $s \le t$ and
	\item $\Op{s}\cdot \Op{t} \subseteq \Op{s+t}$.
	\end{enumerate}
In particular, $\Op{0} \subset \Op{}$ is a subalgebra and $\Op{-\infty} \subset \Op{}$ and $\Op{t} \subset \Op{0}$, $t \in [-\infty,  0)$ are two-sided ideals.
\qed
\end{lem}
Notice that operators with analytic order at most $0$ extend, in particular, to {\em bounded} linear operators on $\H = \W{0}$ allowing us to identify $\Op{0}$ with a subalgebra of $\Linear(\H)$, the algebra of bounded linear operators on $\H$.

\begin{ex} We see from Lemma~\ref{lem: order of Delta^{z}} that the operator $\Delta^{z}$ belongs to $\Op{2\Real(z)}$ for any $z \in \C$. It follows from the spectral theorem that for $z$, $w \in \C$
	\begin{equation}\Delta^{z} \cdot \Delta^{w} = \Delta^{z + w} \quad\text{in } \Op{}.\end{equation}
\end{ex}

\begin{lem}\label{lem: PDOt = Deltahalft PDO0} If $\PDO$ is a filtered subalgebra of $\Op{}$ such that $\Delta^{\half{t}}$ belongs to $\PDO^{t}$ for all $t \in \R$, then $\PDO^{0}$ is unital and
	\begin{equation}\PDO^{t} = \Delta^{\half{t}}\PDO^{0} = \PDO^{0}\Delta^{\half{t}}.\end{equation}
Conversely, if a unital subalgebra $\PDO^{0} \subseteq \Op{0}$ satisfies $\Delta^{\half{t}}\PDO^{0} = \PDO^{0}\Delta^{\half{t}}$ for all $t \in \R$, then $\PDO^{t} \coloneqq \Delta^{\half{t}}\PDO^{0}$ defines a filtered subalgebra of $\Op{}$ such that $\Delta^{\half{t}}$ belongs to $\PDO^{t}$. 
\end{lem}
\begin{proof}
For the first statement: since $\Delta^{0} = 1$, $\PDO^{0}$ is unital. Moreover, for any $t \in \R$,
	\begin{equation}\PDO^{t} = \Delta^{\half{t}}\Delta^{-\half{t}}\PDO^{t} \subseteq \Delta^{\half{t}}\PDO^{-t}\PDO^{t} \subseteq \Delta^{\half{t}}\PDO^{0} \subseteq \PDO^{t}\PDO^{0} \subseteq \PDO^{t}.\end{equation}
Similarly for the other side. 

For the second statement: since $\PDO^{0}$ is unital,
	\begin{align}
	&\PDO^{t} \PDO^{s} = \Delta^{\half{t}}\PDO^{0} \Delta^{\half{s}}\PDO^{0} = \Delta^{\half{t}}\Delta^{\half{s}}\PDO^{0}\PDO^{0} = \PDO^{t+s} \quad \text{and}\\
	&\Delta^{\half{t}} ~\text{belongs~ to}~  \PDO^{t}.
	\end{align}
\end{proof}

\begin{cor}\label{lem: Op^{t}}For any $t \in \R$, 
	\begin{equation}\Op{t} = \Delta^{\half{t}}\Op{0} = \Op{0}\Delta^{\half{t}}.\end{equation}
\qed
\end{cor}

\begin{cor}\label{cor: summability of Op} Let $\Delta$ be an invertible positive self-adjoint operator. Let $\Compact$, $\Linear^{p}$ and $\Linear^{(p, \infty)}$ denote the ideal of compact, $p$-Schatten and $p$-Dixmier operators.
\begin{enumerate}[(a)] 
\item If $\Delta^{-\half{1}} \in \Compact$ then $\Op{-t} \subseteq \Compact$ for $t > 0$.
\item If $\Delta^{-\half{1}} \in \Linear^{p},\: p \ge 1$ then $\Op{-t} \subseteq \Linear^{p/t}$ for $0 < t \le p$.
\item If $\Delta^{-\half{1}} \in \Linear^{(p, \infty)}, \: p \ge 1$ then $\Op{-t} \subseteq \Linear^{(p/t, \infty)}$ for $0 < t \le p$.
\end{enumerate}
\end{cor}
\begin{proof} Using the fact that $\Op{-t} \subseteq \Delta^{-\half{t}}\Linear$, these follow immediatly from well-known properties of the respective ideals.
\end{proof}

\section{Algebra of Generalized Differential Operators}\label{Algebra of Generalized Differential Operators}

In this section, we recall Higson's notion of an algebra of generalized differential operators \cite{H:2004a, H:2006}. Let $\N = \Z_{\ge 0}$.

\begin{defn}[Higson] An {\em $\N$-filtered} subalgebra $\DO \subseteq \Op{}(\Delta)$ is called an {\em algebra of generalized differential operators}\index{algebra!of generalized differential operators} if $\DO$ is closed under the derivation $[\Delta, -]$ and satisfies
	\begin{equation}
	[\Delta, \DO^{k}] \subseteq \DO^{k+1}, \quad k \in \N.
	\end{equation}
\end{defn}

\begin{ex}\label{ex cptly supp diff op} Let $M$ be a complete Riemannian manifold and let $S \to M$ be a Hermitian vector bundle. Let $\Delta$ be the closure of a scalar Laplacian $+1$ (to ensure invertibility). 
Then the algebra $\DO = \DO_c(M, S)$ of compactly supported differential operators acting on the sections of $S$ is an example of an algebra of generalized differential operators. Note that in this example $\Delta$ is not an element of $\DO$. See \cite{HR:2000a}.
\end{ex}

\begin{ex}[Polynomial Weyl Algebra]\label{Weyl algebra}\index{Weyl algebra} Consider the usual Lebesgue measure on $\R^{n}$ and let $\H \coloneqq L^{2}(\R^{n})$. The operator 
	\begin{equation}1 + \sum_{i=1}^{n}\left(x_{i}^{2} - \frac{\partial^{2}}{\partial x_{i}^{2}}\right)\end{equation}
with domain $C^{\infty}_{c}(\R^{n})$ of compactly supported smooth functions is called the {\em harmonic oscillator}\index{operator!harmonic oscillator}. It is essentially self-adjoint 
and strictly positive.
 
Let $\Delta$ denote the closure. Then $\Delta$ is invertible positive self-adjoint operator with compact resolvent. 
Let $\mathcal{W} = \mathcal{W}_{n}$ be the algebra of {\em polynomial} differential operators on $\R^{n}$ acting on $L^{2}(\R^{n})$. As an algebra, it is generated by $x_{1}, \dots, x_{n}$ and $\frac{\partial}{\partial x_{1}}, \dots, \frac{\partial}{\partial x_{n}}$. We filter $\mathcal{W}$ by requiring that
	\begin{equation}\ord(x_{i}) \coloneqq 1 \quad\text{and}\quad\ord\left(\frac{\partial}{\partial x_{i}}\right)\coloneqq 1, \quad 1 \le i \le n.\end{equation}
(The nonzero degree of $x_{i}$ compensates the noncompactness of $\R^{n}$.) Then $\mathcal{W}$ is  an algebra of generalized differential operators and $\Delta$ is an element of $\mathcal{W}^{2}$. See \cite{Shu:2001}.
\end{ex}

Now we relate generalized differential operators to regularity. Let $\AHD$ be a spectral triple and let $\Delta = \D^{2} + 1$. Then it is clear that $\D \in \Op{1}(\Delta) \subset \End(\W{\infty})$. Suppose that $\A \cdot \W{\infty} \subseteq \W{\infty}$. 
Define an $\N$-filtered algebra $\DO \subseteq \End(\W{\infty})$ inductively by
\begin{enumerate}
\item $\DO^{0} \coloneqq $ the subalgebra generated by $\A + [\D, \A] \subseteq \End(\W{\infty})$ and
\item $\DO^{1} \coloneqq \DO^{0} + [\Delta, \DO^{0}] + \DO^{0}[\Delta, \DO^{0}] \subseteq \End(\W{\infty})$ and
\item $\DO^{k} \coloneqq \DO^{k-1} + \sum_{j=1}^{k-1} \DO^{j}\cdot \DO^{k-j} + [\Delta, \DO^{k-1}] + \DO^{0}[\Delta, \DO^{k-1}] \subseteq \End(\W{\infty})$ for $k \ge 2$. 
\end{enumerate}

\begin{thm}[{Higson \cite[Theorem 4.26]{H:2006}}]\label{Higson Regularity} Let $\AHD$ be a spectral triple and let $\Delta = \D^{2} + 1$. Suppose that $\A \cdot \W{\infty} \subseteq \W{\infty}$. Then $\AHD$ is regular if and only if $\DO^{k} \subseteq \Op{k}$, $k \ge 0$.
\end{thm}

We will give a proof in Section~\ref{section Regularity of ST}, essentially rephrasing the arguments in \cite{CM:1995, H:2006} in the language of algebra of generalized pseudo-differential operators, which we describe next.

\section{Algebra of Generalized Pseudodifferential Operators}\label{Algebra of Generalized Pseudodifferential Operators}

In this section, we define and study algebras of generalized pseudodofferential operators. Just as in the classical case, these provide a convenient framework to study index theoretic problems, see \cite{H:2006}.

\subsection{Operators of all orders}

\begin{defn} 
An {\em $\R$-filtered} subalgebra $\PDO \subseteq \Op{}(\Delta)$ is called an {\em algebra of generalized pseudo-differential operators} if $\PDO$ satisfies for $z \in \C$, $t \in \R$,
	\begin{equation}\label{prod Deltas}
	\Delta^{\half{z}}\PDO^{t}\subseteq \PDO^{\Real(z) + t} \quad\text{and}\quad \PDO^{t}\Delta^{\half{z}} \subseteq \PDO^{\Real(z) + t}
	\end{equation}
and
	\begin{equation}\label{comm Deltas}
	[\Delta^{\half{z}}, \PDO^{t}] \subseteq \PDO^{\Real(z) + t - 1}. 
	\end{equation}
\end{defn}

\begin{rem}\label{rem Deltazhalf belong to pdorealz} If $\Delta^{\half{z}}$ belongs to $\PDO^{\Real(z)}$ for $z \in \C$, then (\ref{prod Deltas}) is trivially satisfied and
	\begin{equation}
	\PDO^{t} = \Delta^{\half{t}}\PDO^{0} = \PDO^{0}\Delta^{\half{t}}, \quad t \in \R,
	\end{equation}
by Lemma~\ref{lem: PDOt = Deltahalft PDO0}.
\end{rem}
The following is the classical example.
\begin{ex}[Pseudodifferential Operators]\label{ex: pdo's} Let $M$ be a closed manifold.
Let $\PDO = \bigcup \PDO^{t}$ denote the $\R$-filtered algebra of {\em pseudo-differential operators} on $M$ with scalar symbols and let $\CL \subset \PDO$ denote the filtered subalgebra of {\em classical pseudo-differential operators} on $M$ (see for instance \cite{Shu:2001}).

Let $\Delta$ be (the closure of) an invertible, positive, order two, partial-differential operator on $M$. Then we have inclusions of {\em $\R$-filtered} algebras
	\begin{equation}
	\CL \subset \PDO \subset \Op{}(\Delta).
	\end{equation}
Moreover $\Delta^{\half{z}}$ belongs to $\CL^{\Real(z)} \subset \PDO^{\Real(z)}$, $z \in \C$, and since
	\begin{equation}\label{commutator of pdo}
	[\PDO^{s}, \PDO^{t}] \subseteq \PDO^{s+t-1},\quad s, t \in \R,
	\end{equation}
we have 
	\begin{align}
	\label{diff cl}[\Delta^{\half{z}}, \CL^{t}] \subseteq \CL^{\Real(z) + t - 1}\quad&\text{and}\quad
	[\Delta^{\half{z}}, \PDO^{t}] \subseteq \PDO^{\Real(z) + t -1}.
	\end{align}
Hence $\CL$ and $\PDO$ are algebras of generalized pseudo-differential operators (in view of Remark~\ref{rem Deltazhalf belong to pdorealz}). 

More generally, we could consider pseudo-differential operators acting on sections of a vector bundle. Then we need to assume that $\Delta$ has scalar principal symbol: the property (\ref{commutator of pdo}) would not hold in general, but (\ref{diff cl}) would still hold.

For a calculus of pseudo-differential and classical pseudo-differential operators on $\R^{n}$, see \cite[Chapter IV]{Shu:2001} or \cite[Chapter 3]{Fed:1996}. 
\end{ex}

\begin{rem}\label{rem: DO from PDO} Let $\PDO \subseteq \Op{}$ be an algebra of generalized pseudo-differential operators. Then $\DO^{k} \coloneqq \PDO^{k}$, $k \in \N$, is clearly an algebra of generalized {\em differential} operators. 
\end{rem}
In the other direction, we have the following lemma.
\begin{lem}[{cf.\ \cite[Proposition 4.31]{H:2006}}]\label{lem: PDO from DO} Let $\DO$ be an algebra of generalized differential operators. Let $\PDO^{t}$ denote the space of linear combinations of operators $P \in \Op{}$ such that for any $l \in \R$, $P$ may be decomposed as
	\begin{equation}P = X\Delta^{\half{z-m}} + Q\end{equation}
with $\Real(z) \le t$ and $X \in \DO^{m}$, $m \in \N$ and  $Q \in \Op{l}$. Then $\PDO$ is an algebra of generalized pseudo-differential operators.
\end{lem}

	
First we prove an auxiliary lemma (cf. \cite[Theorem B1]{CM:1995}, \cite[Lemma 4.20]{H:2006}).
	
\begin{notation} In a filtered space, we write
	\begin{equation}P \sim \sum_{k \in \Lambda} P_{k}\end{equation}
if for any $l \in \R$, there exists a finite subset $F \subset \Lambda$ such that $P -\sum_{k \in F} P_{k}$ has order at most $l$.
\end{notation}

\begin{lem}[{cf.\ \cite[Theorem B1]{CM:1995}, \cite[Lemma 4.20]{H:2006}}]\label{lem: Taylor} Let $z \in \C$ and let $Y \in \DO$. Let $Y^{(0)} \coloneqq Y$ and $Y^{(k)} \coloneqq [\Delta, Y^{(k-1)}]$, $k \ge 1$. Then
	\begin{equation}\Delta^{z}Y \sim \sum_{k=0}^{\infty}{z \choose k}Y^{(k)}\Delta^{z - k}.\end{equation} 
\end{lem}
Note that $\ord(Y^{(k)}) \le \ord(Y) + k$, so $\ord (Y^{(k)}\Delta^{z-k}) \le \ord (Y) + 2\Real(z) - k$.
\begin{proof} Assume $\Real(z) < 0$. Let $R$ denote the resolvent $(\lambda - \Delta)^{-1} \in \Op{-2}$. Then it is easy to see that
	\begin{equation}
	RY = YR + RY^{(1)}R.
	\end{equation}
Hence, for any $n \in \N$,
	\begin{equation}RY = YR + Y^{(1)}R^{2} + \dots +Y^{(n)}R^{n+1} + RY^{(n+1)}R^{n+1}.\end{equation}	
Applying the Cauchy integral formula (Lemma~\ref{lem: Cauchy integral formula for complex power}), we see that
	\begin{align}
	\Delta^{z}Y = &\sum_{k=0}^{n}{z \choose k}Y^{(k)}\Delta^{z - k}\\ 
	&+ \frac{1}{2\pi i}\int \lambda^{z} (\lambda - \Delta)^{-1}Y^{(n+1)}(\lambda - \Delta)^{-n-1}d\lambda.\notag
	\end{align}
But the last integral converges absolutely in $\Linear(\W{s+l}, \W{s})$ for $l = \ord(Y) + n + 1 - 2(n+1)$, hence has order at most $\ord(Y) - n - 1$. 

For general $z \in \C$, the identity
	\begin{align*}
	\Delta^{z+1}Y = \Delta^{z} Y \Delta + \Delta^{z}Y^{(1)}
	\end{align*}
allows one to reduce to the case $\Real(z) < 0$.  	
\end{proof}

\begin{proof}[Proof of Lemma~\ref{lem: PDO from DO}] Clearly $\PDO$  is a filtered subspace of $\Op{}$. It follows from Lemma~\ref{lem: Taylor} that $\PDO$ is a subalgebra satisfying (\ref{prod Deltas}). For $P = X\Delta^{\half{w - m}} + Q$ in $\PDO^{t}$ with $\Real(w) = t$ and $X \in \DO^{m}$ and $\ord(Q)$ small, the expansion
	\begin{equation}[\Delta^{\half{z}}, P] \sim \sum_{k=1}^{\infty}{\half{z} \choose k}X^{(k)}\Delta^{\half{z + w - m -2k}} + [\Delta^{\half{z}}, Q]\end{equation}
shows that $[\Delta^{\half{z}}, P] \in \PDO^{\Real(z) + t -1}$, proving (\ref{comm Deltas}).
\end{proof}

\subsection{Operators of order at most zero}
\begin{defn} A subalgebra $\BPDO \subseteq \Op{0}$ is called an {\em algebra of generalized pseudo-differential operators of order at most zero} if $\BPDO$ is closed under the derivation $\delta \coloneqq [\Delta^{\half{1}}, -]$.
\end{defn}

\begin{lem}[{cf.\ \cite[Lemma B1]{CM:1995}}]\label{lem deltainfty} Let $\BPDO = \dom^{\infty}(\delta)$ denote the smooth domain of $\delta$. Then $\BPDO$ is an algebra of generalized pseudo-differential operators of order at most zero.
\end{lem}
\begin{proof}We just need to show that $\BPDO$ is a subset of $\Op{0}$, since by construction $\BPDO$ is an algebra closed under $\delta$.

Let $s = \Real(z) < 0$. Then for any $b \in \BPDO$ and $n \in \N$,
	\begin{align}\label{Deltazb}
	\Delta^{\half{z}}b = \sum_{k=0}^{n}&{z \choose k}\delta^{k}(b)\Delta^{\half{z - k}} \\
	&+ \frac{1}{2\pi i}\int \lambda^{z} (\lambda - \Delta^{\half{1}})^{-1}\delta^{n+1}(b)(\lambda - \Delta^{\half{1}})^{-n-1}d\lambda,\notag
	\end{align}
in $\Linear(\H)$, where the integral is a contour integral along a downwards pointing vertical line in $\C$ which separates $0$ from $\Spec(\Delta^{\half{1}})$. The proof of Lemma~\ref{lem: Taylor} applied to $\Delta^{\half{1}}$ goes through ad verbatim.

The last integral gives a bounded operator in $\Linear(\H, \W{n+1})$ so taking $n$ large enough we see that $||\Delta^{\half{s}}b\xi|| \le C\cdot ||\Delta^{\half{s}}\xi||$, $\xi \in \dom(\Delta^{\half{s}})$, for some $C > 0$. Hence $b \cdot \W{s} \subseteq \W{s}$ for $s < 0$. For $s > 0$, the identity
	\begin{equation}\Delta^{\half{s+1}}b = \Delta^{\half{s}}b\Delta^{\half{1}} + \Delta^{\half{s}}\delta(b)\end{equation}
allows one to reduce to the case $s < 0$. Hence $\BPDO \subseteq \Op{0}$. 

\end{proof}

\begin{rem}\label{BPDO from PDO} It is clear that if $\PDO \subseteq \Op{}$ is an algebra of generalized pseudo-differential operators then $\BPDO \coloneqq \PDO^{0}$ is an algebra of generalized pseudo-differential operators of order at most zero.
\end{rem}

Conversely, we have the following.

\begin{lem}[{cf.\ \cite[Lemma 4.27]{H:2006}}]\label{BPDO to DO} Let $\BPDO \subseteq \Op{0}$ be an algebra of generalized pseudo-differential operators of order at most zero. Then 
	\begin{equation}\DO^{k} \coloneqq \sum_{j=0}^{k}\Delta^{\half{j}}\BPDO \subseteq \Op{k},\quad k \in \N,
	\end{equation}
is an algebra of generalized differential operators.
\end{lem}
\begin{proof}
First note that since 
	\begin{equation}\label{Delta BPDO}
	\Delta^{\half{k}}b = \sum_{j = 0}^{k} {k \choose j}\delta^{j}(b)\Delta^{\half{k-j}}, \quad k \in \N,  
	\end{equation}
we have
	\begin{equation}
	\Delta^{\half{k}}\BPDO = \BPDO\Delta^{\half{k}}, \quad k \in \N.
	\end{equation}
Hence $\DO$ is an $\N$-filtered subalgebra of $\Op{}$ (see the proof of Lemma~\ref{lem: PDOt = Deltahalft PDO0}).
Then using the facts 	
	\begin{align}
	\label{Deltahalf, Deltakb}
	\delta(\Delta^{\half{k}}\BPDO) &= \Delta^{\half{k}}\delta(\BPDO) \subseteq \Delta^{\half{k}}\BPDO\quad \text{and}\\
	 \label{Delta commutator delta}[\Delta, P] &= 
	2\Delta^{\half{1}}\delta(P) - \delta(\delta(P)),
	\end{align}
we see that $\DO$ is indeed an algebra of generalized differential operators.
\end{proof}

\begin{rem} Suppose that $\BPDO$ is an algebra of generalized pseudo-differential operators of order at most zero satisfying, for all $t \in \R$, 
	\begin{equation}
	\Delta^{it} \in \BPDO \quad \text{and} \quad\Delta^{\half{t}}\BPDO\Delta^{-\half{t}} \subseteq \BPDO.
	\end{equation}
Then one can show that $\PDO^{t} \coloneqq \Delta^{\half{t}}\BPDO$, $t \in \R$,  is an algebra of generalized pseudo-differential operators, using Lemma~\ref{lem: Taylor}. This is applicable to $\BPDO = \dom^{\infty}(\delta)$ of Lemma~\ref{lem deltainfty}.
\end{rem}

\section{Regularity of Spectral Triples}\label{section Regularity of ST}
See \cite[Appendix B]{CM:1995}, \cite[Theorem 3.25]{H:2004a}, \cite[Theorem 4.26]{H:2006}.
\begin{thm}\label{thm regularity} Let $\AHD$ be a spectral triple. Let $\Delta \coloneqq \D^{2} + 1$ and let $\delta \coloneqq [\Delta^{\half{1}}, -]$. Then the following conditions are equivalent:
\begin{enumerate}[(1)]
\item\label{reg |D|} The spectral triple $\AHD$ is regular, that is, $\A + [\D, \A]$ is contained in $\dom^{\infty}([|\D|, -])$.
\item\label{reg delta} The set $\A + [\D, \A]$ is contained in $\dom^{\infty}(\delta)$.
\item\label{reg BPDO} There exists an algebra of generalized pseudo-differential operators of order at most zero containing $\A + [\D, \A]$.
\item\label{reg DO} There exists an algebra of generalized differential operators containing $\A + [\D, \A]$ in degree $0$.
\item\label{reg PDO} There exists an algebra of generalized pseudo-differential operators containing $\A + [\D, \A]$ in degree $0$.
\end{enumerate}
\end{thm}
\begin{proof}
(\ref{reg |D|}) $\Leftrightarrow$ (\ref{reg delta}) is clear, since $\Delta^{\half{1}}$ is a bounded perturbation of $|\D|$. The implication (\ref{reg delta}) $\Rightarrow$ (\ref{reg BPDO}) follows from Lemma~\ref{lem deltainfty}.

To prove (\ref{reg BPDO})  $\Rightarrow$ (\ref{reg delta}), suppose that $\BPDO \subseteq \Op{0}$ is an algebra of generalized pseudo-differential operators of order at most zero. Then $\BPDO$ is contained in $\dom^{\infty}(\delta)$. Indeed, if $b \in \BPDO \subseteq \Op{0}$, then $b \cdot \dom(\Delta^{\half{1}}) \subseteq \dom(\Delta^{\half{1}})$ and, being an operator of order $0$, the commutator $[\Delta^{\half{1}}, b] \in \BPDO \subseteq \Op{0}$ extends to a bounded operator on $\H$. In other words, $b$ belongs to $\dom(\delta)$. But $\BPDO$ is closed under $\delta$, therefore $\BPDO \subseteq \dom^{\infty}(\delta)$.

The implication (\ref{reg BPDO}) $\Rightarrow$ (\ref{reg DO}) follows from Lemma~\ref{BPDO to DO}, (\ref{reg DO}) $\Rightarrow$ (\ref{reg PDO}) follows from Lemma \ref{lem: PDO from DO} and (\ref{reg PDO}) $\Rightarrow$ (\ref{reg BPDO}) follows from Remark \ref{BPDO from PDO}. 
\end{proof}

As a corollary we obtain a proof of Theorem~\ref{Higson Regularity}.
\begin{proof}[Proof of Theorem \ref{Higson Regularity}] If $\AHD$ is regular, then by Theorem~\ref{thm regularity}, there exists an algebra of generalized differential operators $\mathcal{E} \subseteq \Op{}$. It is clear, by induction, that $\DO^{k} \subseteq \mathcal{E}^{k}$, $k \in \N$. Thus $\DO^{k} \subseteq \Op{k}$. Conversely, if $\DO^{k} \subseteq \Op{k}$, $k \in \N$, then $\DO$ is an algebra of generalized differential operators and by Theorem~\ref{thm regularity}, $\AHD$ is regular.
\end{proof}

\begin{ex} \begin{enumerate}
\item The commutative spectral triple of Example~\ref{comm example} is regular, because of Example~\ref{ex cptly supp diff op} or Example~\ref{ex: pdo's}.
\item The spectral triple of \cite{CM:1995} associated to a triangular structure is regular because the algebra of $\Psi$DO$'$-operators is an example of an algebra of generalized pseudo-differential operators.
\end{enumerate}
\end{ex}

\begin{lem}\label{prop product} Let $(\A_{1}, \H_{1}, \D_{1})$ and $(\A_{2}, \H_{2}, \D_{2})$ be spectral triples. Let $\A_{1} \otimes_{\alg} A_{2}$ denote the algebraic tensor product and let $\H_{1} \tensor \H_{2}$ denote the graded Hilbert space tensor product. Then the operator 
	\begin{equation}
	\D \coloneqq \D_{1} \tensor 1 + 1 \tensor \D_{2}
	\end{equation}
with domain 
	\begin{equation}
	\dom(\D) \coloneqq \dom(\D_{1}) \tensor_{\alg} \dom(\D_{2}) \subseteq \H,
	\end{equation}
is essentially self-adjoint and  $(\A_{1} \otimes \A_{2}, \H_{1} \tensor \H_{2}, \bar\D)$ is a spectral triple.
\end{lem}
\begin{proof} See \cite{otgonbayar-2009}.
\end{proof}

We write $\D_{1} \times \D_{2}$ for the closure $\bar\D$ of $\D$ and call 
	\begin{equation}\label{defn eq product}
	(\A_{1} \otimes \A_{2}, \H_{1} \tensor \H_{2}, \D_{1} \times \D_{2}).
	\end{equation}
the product of $(\A_{1}, \H_{1}, \D_{1})$ and $(\A_{2}, \H_{2}, \D_{2})$.

The following theorem is implicitly contained in or follows from results in \cite{CM:1995, H:2004a, H:2006}. But no direct reference seems to exist in the literature.

\begin{thm}\label{thm: product of regular} The product of regular spectral triples are again regular. 
\end{thm}
\begin{proof} Let $(\A_{1}, \H_{1}, \D_{1})$ and $(\A_{2}, \H_{2}, \D_{2})$ be regular spectral triples and let $\DO_{i} \subseteq \Op{}(\H_{i}, \D_{i}^{2} + 1)$ be algebras of generalized differential operators for $(\A_{i}, \H_{i}, \D_{i})$ respectively, i.e.\ $\A_{i} + [\D_{i}, \A_{i}] \subseteq \DO_{i}^{0}$ for $i = 1, 2$. 

We claim that $\DO_{1} \tensor \DO_{2}$, with the product filtering, is an algebra of generalized differential operators for the product $(\A_{1} \otimes \A_{2}, \H_{1} \tensor \H_{2}, \D_{1} \times \D_{2})$. Then by Theorem~\ref{thm regularity}, we see that the product is regular. 

Let $\Delta_{i} = \D_{i}^{2} + 1$ and $\Delta = (\D_{1} \times \D_{2})^{2} + 1$. Then it is easy to check that
	\begin{equation}
	\Op{s}(\H_{1}, \Delta_{1}) \tensor \Op{t}(\H_{2}, \Delta_{2}) \subseteq \Op{s+t}(\H_{1} \tensor \H_{2}, \Delta).
	\end{equation}
Hence $\DO_{1} \tensor \DO_{2}$ is a filtered subalgebra of $\Op{}(\H_{1} \tensor \H_{2}, (\D_{1} \times \D_{2})^{2} + 1)$. For $P_{1} \in \DO_{1}^{k}$, $P_{2} \in \DO_{2}^{l}$, we see that
	\begin{equation}
	[\Delta, P_{1} \tensor P_{2}] = [\Delta_{1}, P_{1}] \tensor P_{2} + P_{1} \tensor [\Delta_{2}, P_{2}],
	\end{equation}
belongs to $(\DO_{1} \tensor \DO_{2})^{k + l - 1}$ and thus $\DO_{1} \tensor \DO_{2}$ is an algebra of generalized differential operators. Finally, for $a_{1} \in \A_{1}$ and $a_{2} \in \A_{2}$, clearly $a_{1} \otimes a_{2}$ belongs to $(\DO_{1} \tensor \DO_{2})^{0}$ and so does
	\begin{equation}
	[\D_{1} \times \D_{2}, a_{1} \otimes a_{2}] = [\D_{1}, a_{1}] \tensor a_{2} + a_{1} \tensor [D_{2}, a_{2}]. 
	\end{equation}
This completes the proof.
\end{proof}

\bibliographystyle{amsalpha}
\bibliography{regularity}

\providecommand{\bysame}{\leavevmode\hbox to3em{\hrulefill}\thinspace}
\providecommand{\MR}{\relax\ifhmode\unskip\space\fi MR }
\providecommand{\MRhref}[2]{%
  \href{http://www.ams.org/mathscinet-getitem?mr=#1}{#2}
}
\providecommand{\href}[2]{#2}
\begin{thebibliography}{Shu01}

\bibitem[CM95]{CM:1995}
A.~Connes and H.~Moscovici, \emph{The local index formula in noncommutative
  geometry}, Geom. Funct. Anal. \textbf{5} (1995), no.~2, 174--243.
  \MR{MR1334867 (96e:58149)}

\bibitem[Fed96]{Fed:1996}
Boris Fedosov, \emph{Deformation quantization and index theory}, Mathematical
  Topics, vol.~9, Akademie Verlag, Berlin, 1996. \MR{MR1376365 (97a:58179)}

\bibitem[Hig04]{H:2004a}
Nigel Higson, \emph{The local index formula in noncommutative geometry},
  Contemporary developments in algebraic $K$-theory, ICTP Lect. Notes, XV,
  Abdus Salam Int. Cent. Theoret. Phys., Trieste, 2004, pp.~443--536
  (electronic). \MR{MR2175637}

\bibitem[Hig06]{H:2006}
\bysame, \emph{The residue index theorem of {C}onnes and {M}oscovici}, Surveys
  in noncommutative geometry, Clay Math. Proc., vol.~6, Amer. Math. Soc.,
  Providence, RI, 2006, pp.~71--126. \MR{MR2277669}

\bibitem[HR00]{HR:2000a}
Nigel Higson and John Roe, \emph{Analytic {$K$}-homology}, Oxford Mathematical
  Monographs, Oxford University Press, Oxford, 2000, Oxford Science
  Publications. \MR{MR1817560 (2002c:58036)}

\bibitem[Otg09]{otgonbayar-2009}
Uuye Otgonbayar, \emph{Multiplicativity of the {JLO}-character}, 2009.

\bibitem[RS75]{RS:1975}
Michael Reed and Barry Simon, \emph{Methods of modern mathematical physics.
  {II}. {F}ourier analysis, self-adjointness}, Academic Press [Harcourt Brace
  Jovanovich Publishers], New York, 1975. \MR{MR0493420 (58 \#12429b)}

\bibitem[Shu01]{Shu:2001}
M.~A. Shubin, \emph{Pseudodifferential operators and spectral theory}, second
  ed., Springer-Verlag, Berlin, 2001, Translated from the 1978 Russian original
  by Stig I. Andersson. \MR{MR1852334 (2002d:47073)}

\end{thebibliography}
\end{document}